\documentclass[11pt,onecolumn,journal]{IEEEtran}

\usepackage[cmex10]{amsmath}
\usepackage{amssymb,amsthm}
\usepackage{graphicx}
\usepackage{bm,xspace,fancyhdr}
\usepackage[tight,footnotesize]{subfigure}
\usepackage{epstopdf}
\usepackage{verbatim}
\usepackage{cite,url}
\usepackage{pdfsync}
\usepackage{setspace}

\newtheorem{lemma}{Lemma}[section]
\newtheorem{thm}{Theorem}[section]

\newcommand{\bmath}[1]	{\ensuremath{\bm{#1}}\xspace}

\newcommand{\xmath}[1]	{\ensuremath{#1}\xspace}

\newcommand{\reals}	{\xmath{\mathbb{R}}}
\newcommand{\norm}[1]	{\xmath{\left|\left| {#1} \right|\right|}}

\newcommand{\lz}	{\xmath{\ell_0}}
\newcommand{\lo}	{\xmath{\ell_1}}

\newcommand{\stepsz}	{\xmath{\mu}}

\newcommand{\fvecmat}	{\xmath{\bm{\Phi}}}

\newcommand{\basemat}	{\xmath{\bm{W}}}
\newcommand{\ripd}	{\xmath{\delta}}

\newcommand{\ripc}	{\xmath{C}}

\newcommand{\istavar}	{{u}}
\newcommand{\istavec}	{\xmath{\bm{\istavar}}}
\newcommand{\istavect}[1]	{\xmath{\bm{\istavar}_{#1}}}

\newcommand{\coefsym}	{{z}}
\newcommand{\coefvec}	{\xmath{\bm{\coefsym}}}
\newcommand{\coefprev}	{\xmath{\widetilde{\bm{\coefsym}}}}
\newcommand{\coeftrue}	{\xmath{\bm{\coefsym}^{\dagger}}}

\newcommand{\coefvect}[1]	{\xmath{\bm{\coefsym}_{#1}}}
\newcommand{\coefvest}	{\xmath{\widehat{\bm{\coefsym}}}}
\newcommand{\coefvtest}[1]	{\xmath{\widehat{\bm{\coefsym}}_{#1}}}

\newcommand{\insym}	{y}
\newcommand{\insig}	{\xmath{\bm{\insym}}}
\newcommand{\insigt}[1]	{\xmath{\bm{\insym}_{#1}}}

\newcommand{\measdim}	{\xmath{M}}
\newcommand{\statedim}	{\xmath{N}}

\newcommand{\sparsity}	{\xmath{S}}

\newcommand{\statesym}	{x}
\newcommand{\statev}	{\xmath{\bm{\statesym}}}
\newcommand{\statevt}[1]	{\xmath{\bm{\statesym}_{#1}}}
\newcommand{\statevtest}[1]	{\xmath{\widehat{\bm{\statesym}}_{#1}}}

\newcommand{\statevest}	{\xmath{\widehat{\bmath{\statesym}}}}
\newcommand{\statevestt}[1]	{\xmath{\widehat{\bmath{\statesym}}_{#1}}}

\newcommand{\dfunsym}	{f}

\newcommand{\dfunct}[2]	{\xmath{\dfunsym_{#2}\left(#1\right)}}
\newcommand{\dlip}	{\xmath{\dfunsym^{\ast}}}

\newcommand{\memat}	{\xmath{\bm{G}}}
\newcommand{\mematt}[1]	{\xmath{\bm{G}_{#1}}}

\newcommand{\innovsym}	{\nu}
\newcommand{\innovv}	{\xmath{\bm{\innovsym}}}
\newcommand{\innovvt}[1]	{\xmath{\bm{\innovsym}_{#1}}}
\newcommand{\inmax}	{\xmath{\overline{\innovsym}}}
\newcommand{\meerrsym}	{\epsilon}
\newcommand{\meerrv}	{\xmath{\bm{\epsilon}}}
\newcommand{\meerrvt}[1]	{\xmath{\bm{\epsilon}_{#1}}}
\newcommand{\errmax}	{\xmath{\overline{\meerrsym}}}

\newcommand{\thresh}	{\xmath{\lambda}}
\newcommand{\kthresh}	{\xmath{\kappa}}
\newcommand{\ssize}	{\xmath{\eta}}
\newcommand{\ssizek}	{\xmath{\widetilde{\ssize}}}
\newcommand{\threshvec}	{\xmath{\bm{\lambda}}}

\newcommand{\errvt}[1]	{\xmath{\bm{e}_{#1}}}

\providecommand{\defstart}{
\begin{tabular}{l|l|l} Macro & Symbol & Meaning \\ \hline }
\providecommand{\defstop}{\end{tabular}}

\begin{document}

\doublespacing

\title{Dynamic Filtering of Time-Varying Sparse Signals via \lo Minimization}%
\author{Adam~S.~Charles,~\IEEEmembership{Student Member,~IEEE}
         and~Aurele~Balavoine,~\IEEEmembership{Student Member,~IEEE}
         and~Christopher~J.~Rozell,~\IEEEmembership{Senior Member,~IEEE}
\thanks{Manuscript received June 29, 2015.  This work was supported in part by NSF grant CCF-1409422 and the James S. McDonnell Foundation.}%
\thanks{A. Charles, A. Balavoine and C. Rozell are with the School of Electrical and Computer Engineering, Georgia Institute of Technology, Atlanta, GA, 30332-0250 USA (e-mail: {acharles6,aurele.balavoine,crozell}@gatech.edu). The authors would like to thank S. Asif and A. Kressner for valuable feedback, and C. Cadieu for the database of documentary videos. }}
	   
\maketitle

\begin{abstract}


Despite the importance of sparsity signal models and the increasing prevalence of high-dimensional streaming data, there are relatively few algorithms for dynamic filtering of time-varying sparse signals.  Of the existing algorithms, fewer still provide strong performance guarantees.  This paper examines two algorithms for dynamic filtering of sparse signals that are based on efficient \lo optimization methods. We first present an analysis for one simple algorithm (BPDN-DF) that works well when the system dynamics are known exactly.  We then introduce a novel second algorithm (RWL1-DF) that is more computationally complex than BPDN-DF but performs better in practice, especially in the case where the system dynamics model is inaccurate.  Robustness to model inaccuracy is achieved by using a hierarchical probabilistic data model and propagating higher-order statistics from the previous estimate (akin to Kalman filtering) in the sparse inference process.  We demonstrate the properties of these algorithms on both simulated data as well as natural video sequences.  Taken together, the algorithms presented in this paper represent the first strong performance analysis of dynamic filtering algorithms for time-varying sparse signals as well as state-of-the-art performance in this emerging application.

\end{abstract}



\section{Introduction}
\label{sec:intro}

Many signal processing applications involve performing inference on high-dimensional time-varying signals.  
The past decade of work in signal and image processing has demonstrated that signal priors based on the notion of sparsity (e.g., compressibility in a basis) lead to state-of-the-art performance in many linear inverse problems  (e.g., undersampling, inpainting, denoising, compressed sensing, etc.~\cite{ELA:2008,TAO:2006,BAR:2007}) across several different application domains (e.g. natural images~\cite{ELA:2008}, audio~\cite{DAV:2006}, hyperspectral imagery~\cite{CHA:2011h}, etc.).  In this paradigm, a static signal \statev is observed through a linear measurement process
$\insig = \fvecmat\statev + \meerrv$, and the signal is estimated through performing penalized least-squares with a sparsity-inducing regularization function.  

With the recent increase in streaming data sources, it is critical to develop algorithms that can operate efficiently on these data streams.  In \emph{dynamic filtering}, a causal algorithm estimates the current system state recursively based on the previous state estimate (as opposed to batch methods which jointly estimate the time-varying states once all data is collected).  Classic dynamic filtering approaches stem from Kalman filtering~\cite{KAL:1960}, which has restrictive assumptions (linear systems and Gaussian distributions) and involves large matrix inversion that is inappropriate for high-dimensional data. Given the historical success of dynamic filtering and recent achievements of sparsity models, there are many applications where one can reasonably expect significant gains from jointly leveraging sparsity and models of the signal dynamics.  Unfortunately, compared to the individual literatures on sparse signal estimation and dynamic filtering, there is comparatively little work incorporating both types of information jointly. In particular, existing work (reviewed in detail in Section~\ref{sec:SKFback}) does not integrate the powerful \lo optimization methods from the sparsity literature with the successful Kalman philosophy of using higher-order statistics to propagate information about past states to the current estimate.  Furthermore, existing reports lack strong convergence and performance guarantees.

In this paper we focus on the problem of estimating a time-varying state $\statevt{n}\in\reals^N$ that evolves according to a Markov process
\begin{gather}
	\statevt{n} = \dfunct{\statevt{n-1}}{n} + \innovvt{n},		\label{eq:Process}
\end{gather}
where $\statevt{n}$ obeys a sparsity model at each time step and \dfunct{\cdot}{n}$:\reals^N \rightarrow \reals^N$ is the operator (assumed known) describing the system evolution with some error \innovvt{n}$\in\reals^N$ (called the \emph{innovation}).  The desired state is hidden and only observed through noisy linear measurements $\insig\in\reals^M$:
\begin{gather}
	\insigt{n} = \mematt{n}\statevt{n} + \meerrvt{n}, 		\label{eq:measure}
\end{gather}
where $\mematt{n}\in\reals^{M\times N}$ is the (known but potentially time-varying) measurement matrix\footnote{For static problems we will use \fvecmat as the observation matrix and for dynamic problems we will use \memat in order to better differentiate the two cases.} that has error \meerrvt{n}$\in\reals^M$ and may have $M\ll N$.   

In this paper we examine two methods for \lo-based dynamic filtering of sparse signals. 
We first consider the case where the system dynamics are known with high accuracy. For this case we present a simple and efficient algorithm (BPDN-DF) that combines dynamic penalization with \lo optimization and we present analytic accuracy guarantees on its performance.  
We then introduce a novel second algorithm (RWL1-DF) that performs better in practice at the cost of a marginally higher computational complexity than BPDN-DF. 
In particular, RWL1-DF continues to perform well when the system dynamics model is highly inaccurate.
RWL1-DF achieves robustness to model inaccuracy via a hierarchical probabilistic data model and propagation of higher-order statistics from the previous estimate in the sparse inference process, much as in traditional Kalman filtering.  
We demonstrate the properties of these algorithms both in extensive simulations as well as on natural video sequences.  
In total, we present here both the first strong performance analysis of causal dynamic filtering for temporally evolving sparse signals, as well as state-of-the-art performance in this emerging application.




\section{Background and Related Work} \label{sec:DSest}
\label{sec:background}

\subsection{Sparse Signal Estimation} \label{sec:CSest}
\label{sec:sparsebackground}

In the sparsity model, we assume the signal of interest can be written $\statev = \basemat\coefvec$,
where the coefficient vector $\coefvec\in\reals^{N'}$ is mostly composed of zeros and the dictionary \basemat may be orthonormal or overcomplete.
In the Bayesian interpretation, this model assumes a prior distribution on \coefvec which has \emph{high kurtosis} to encourage coefficients to be zero (or close to zero). Most commonly, an independent, identically distributed (IID)  Laplacian distribution is used as the coefficient prior
$p\left(\bm{z}\right) = \frac{\gamma}{2} e^{-\gamma\norm{\coefvec}_1},$
where $\gamma$ is related to the variance of the individual coefficients and $\norm{\coefvec}_1 = \sum_i\left|\coefvec[i] \right|$ indicates the \lo norm.  Under a Gaussian noise assumption on the measurements $\insig=\fvecmat\statev+ \meerrv$, the corresponding MAP estimate of the signal reduces to the Basis Pursuit De-Noising (BPDN) optimization program
\begin{gather}
	\coefvest = \arg\min_{\coefvec}\left[\frac{1}{2\sigma_{\meerrsym}^2}\norm{\insig - \fvecmat\basemat\coefvec}_2^2 + \gamma\norm{\coefvec}_1 \right], \label{eq:l1optreg} 
\end{gather}
where $\sigma_{\epsilon}^2$ is the  measurement noise variance and the signal estimate is $\statevest = \basemat\coefvest$.
BPDN has been a particularly popular approach due to its strong performance guarantees~\cite{DON:2005} and the development of many specialized optimization approaches~\cite{CAN:2005, BOY:2007,MAL:2005,NOW:2007,ROZ:2008,shapero2012jetcas}.  
As a guarantee of the measurement quality, we say that \fvecmat satisfies the restricted isometry property with parameters 2\sparsity and $\ripd$ (RIP($2\sparsity$,$\ripd$)) with respect to \basemat if for every 2\sparsity-sparse vector \coefvec we have that
 \begin{gather}
 	\ripc(1-\ripd)\norm{\coefvec}_2^2 \leq \norm{\fvecmat\basemat\coefvec}_2^2 \leq \ripc(1+\ripd)\norm{\coefvec}_2^2.  
 \end{gather}
If \coefvec has only \sparsity non-zeros and \fvecmat satisfies RIP($2\sparsity$,$\ripd$) for an appropriate value of $\ripd$, 
then the BPDN estimation error can be shown to be bounded as 
\begin{gather}
	\norm{\coefvec - \coefvest}_2 \leq C_1\norm{\meerrv}_2 + \frac{C_2}{\sqrt{\sparsity}}\norm{\coefvec - \widetilde{\coefvec}_\sparsity}_1, 
\end{gather}
where $C_1$ and $C_2$ are constants and $\widetilde{\coefvec}_\sparsity$ is the best \sparsity-term approximation to \coefvec~\cite{CAN:2007,zhang2009some,aureleISTA}. The RIP condition serves to quantify the quality of the measurements by bounding how similar measurements of different signals could be. 
The choice of basis \basemat is particularly important in ensuring the RIP holds. For orthonormal $\bm{W}$, simply ensuring that the dot products between the rows of \fvecmat and the columns of \basemat are small (i.e. incoherent) is sufficient. 
In the case of overcomplete dictionaries, the inner products between the columns of \basemat must also be small, a property that can still be found in tight frames. 

In~\eqref{eq:l1optreg}, the variance parameter $\gamma$ is assumed known and identical for each coefficient. The RWL1 algorithm~\cite{WAK:2008}  improves performance by generalizing these assumptions through
 a hierarchical probabilistic model known as the Laplacian scale mixture (LSM)~\cite{GAR:2010}.  The LSM introduces a second layer of random variables $\threshvec\in\reals^{N'}$ representing the variance of the first level Laplacian coefficients $\coefvec$,  allowing a variable SNR on each coefficient that is also inferred from the data.  The LSM uses the Gamma distribution (a conjugate prior for the Laplacian) as the hyperprior on $\threshvec[i]$
\begin{gather} 
	p(\threshvec[i]) = \frac{\threshvec^{\alpha-1}[i]}{\theta^{\alpha}\Gamma\left(\alpha\right)}e^{-\threshvec[i]/\theta}, \nonumber
\end{gather}
where $\Gamma\left(\cdot\right)$ is the Gamma function and $\alpha$ and $\theta$ are hyper-parameters with $E[\threshvec[i]] = \alpha\theta$ and $\mbox{Var}[\threshvec[i]] = \alpha\theta^2$. 

Assuming $\alpha = 1$ and a Laplacian conditional distribution on the coefficients $p\left(\bm{z}[i] {|} \threshvec[i]\right) =  \frac{\threshvec[i]}{2} e^{-\threshvec[i]|\bm{z} [i]|}$, the MAP estimate for the full model can be written as 
\begin{gather}
	\left\{\coefvest, \widehat{\threshvec}\right\} = \arg\min_{\coefvec, \threshvec} \norm{\insig - \fvecmat\basemat\coefvec}_2^2 + \sum_i{\left| \threshvec[i]\coefvec[i]\right|}  \nonumber \\
	\qquad\qquad\qquad - \log\left(\sum_i{\threshvec[i]}\right) + \theta^{-1}\sum_i{\threshvec[i]}. \nonumber
\end{gather}
Though this optimization program is difficult to optimize directly, employing expectation-maximization (EM)~\cite{GAR:2010} 
results in the RWL1 algorithm~\cite{WAK:2008} where the following equations are iterated until a convergence criteria is met:
\begin{eqnarray}
	\coefvest^t & = & \arg\min_{\coefvec} \norm{\insig - \fvecmat\basemat\coefvec}_2^2 + \thresh_0\sum_i{\left| \threshvec^t[i]\coefvec[i]\right|} \nonumber \\
        \threshvec^{t+1}[k] & = & \frac{\beta}{|\coefvest^t[k]| + \eta}, \label{eqn:RWL1reg} \nonumber
\end{eqnarray}
with $t$ indicating the iteration number and $\{\thresh_0, \beta, \eta\}$ constant hyperprior parameters. Several papers have followed up on this RWL1 approach, including a description of it as a compromise between \lo and \lz minimization~\cite{WIP:2010,YIN:2008} and an analysis of performance guarantees~\cite{needell2009noisy,KHA:2010,HAS:2010b}. Intuitively, coefficients with little current evidence will have their variances decreased so that future iterations are more likely to estimate them close to zero, and coefficients with strong current evidence will have their variances increased to make future iterations more likely to use those coefficients.  

\subsection{Dynamic Signal Estimation}

\begin{figure*}
	\centering
	\includegraphics[width=\textwidth]{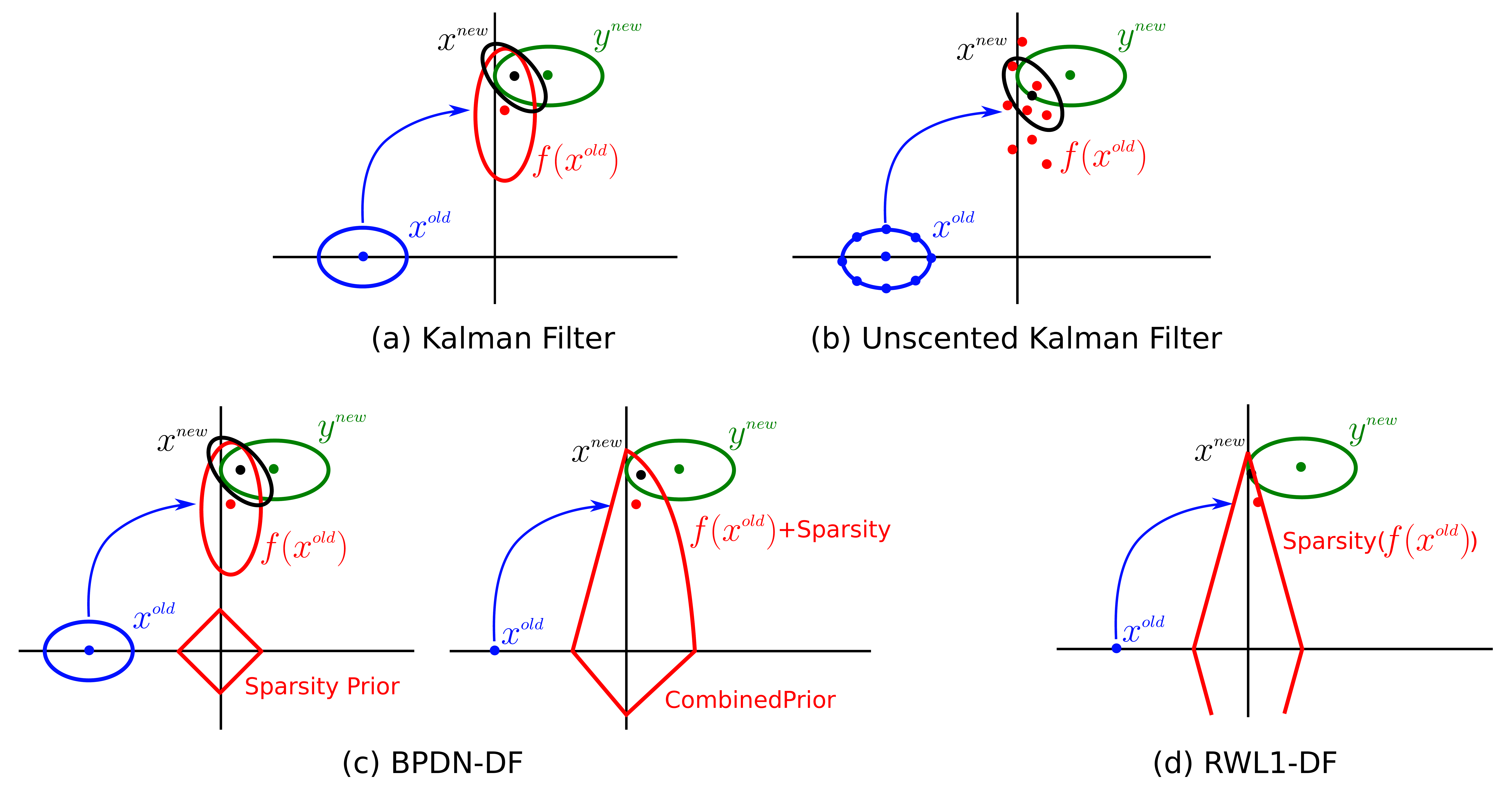}
	\caption{Information propagation in dynamic filtering algorithms.  (a) Classic Kalman filtering propagates the current estimated mean and covariance (blue) through the system dynamics to generate a prior distribution (red) for the next state estimate. This prior is used in conjunction with the likelihood function (green) from the measurements to form a posterior (black) used in estimation. (b) The unscented Kalman filter (and particle filtering generally) estimate the prior distribution via a sampling process that empirically approximates the distribution. (c) BPDN-DF adds a sparsity inducing term to the prior, resulting in a combined prior that is curved outward reminiscent of an $\ell_2$ ball.
	 (d) RWL1-DF uses the previous estimate to set the parameters of a prior that has the diamond-like shape known to promote sparse solutions.
		}
	\label{fig:Style1}
\end{figure*}

In dynamic filtering, knowledge of the system dynamics is used to set a prior probability distribution on the current estimate.
In the classic Kalman filter setup, the linear and Gaussian assumptions (including the system dynamics $\dfunct{\statev}{k} = \bm{F}_k\statev$) translates to a Gaussian assumption on the innovations $\innovvt{k} \sim \mathcal{N}\left(\bm{0}, \bm{Q}\right)$ and on the distribution of the current state conditioned on the previous state
$p(\statevt{k} {|} \statevt{k-1}) \sim \mathcal{N}\left(\bm{F}_k\statevt{k-1}, \bm{Q} \right)$.  Similarly, assuming $\meerrvt{k} \sim \mathcal{N}\left(\bm{0}, \bm{R}\right)$, the likelihood on the measurements at time $k$ conditioned on the current state is
$p(\insigt{k} {|} \statevt{k}) \sim \mathcal{N}\left(\mematt{k}\statevt{k}, \bm{R} \right)$.  Using the Markov property of the state distributions we can write the current state's marginal distribution as $p(\statevt{n}) \sim \mathcal{N}\left(\statevestt{n-1}, \bm{F}_n\bm{P}_{n-1}\bm{F}_n^T + \bm{Q} \right)$, where $\bm{P}_{n}$ is the covariance matrix of the previous state estimate.  The Kalman filter infers the current state via a MAP estimate
\begin{eqnarray}
	\statevtest{n} & = & \arg\min_{\statevt{n}}\left[ \norm{\insigt{n} - \mematt{n}\statevt{n}}_{\bm{R}, 2}^2\right. \nonumber \\
	& & \qquad \left. + \norm{\statevt{n} - \bm{F}_n\statevtest{n-1}}_{\bm{F}_n\bm{P}_{n-1}\bm{F}_n^T + \bm{Q}, 2}^2 \right],  	\label{eq:KalOpt}
\end{eqnarray}
which has an analytic solution (known as the Kalman update equations~\cite{KAL:1960}).  The causal estimation procedure is depicted in Figure~\ref{fig:Style1}(a).  
Intuitively, the estimator propagates the previous state estimate (and its covariance matrix) through the system dynamics to get a prediction (i.e., a prior) on the current state.  This prior is combined with the likelihood resulting from the current measurement to generate a posterior distribution, which is used to estimate the current state. Though the causal estimator is computationally simpler than joint estimation of the states, it still calculates the same estimate for \statevtest{n} as if all of the previous data had been used due to the information propagation at each iteration via the covariance matrices.

Unfortunately, the analytic simplicity and optimality guarantees of the Kalman filter are highly dependent on the linear and Gaussian model assumptions.  
Though not optimal estimators, many heuristic approaches have been introduced that follow the spirit of the Kalman filter while incorporating nonlinear system dynamics (e.g., Extended Kalman Filter~\cite{HAY:2001}) or non-Gaussian structure.   An alternate approach for highly non-linear functions or non-Gaussian statistics is particle filtering, which uses discrete points (particles) to approximate the relevant distributions and to propagate those distributions through nonlinear dynamics functions. The Unscented Kalman filter~\cite{WAN:2001} (depicted in Figure~\ref{fig:Style1}(b)) is similar but uses deterministic (rather than Monte-Carlo) particles. Though particle filtering approaches are general, these methods can become intractable in high-dimensional state spaces and do not explicitly utilize information about sparsity structure. 

\subsection{Sparsity in Dynamic Signals}
\label{sec:SKFback}

Recent work has begun to address time-varying sparse signal estimation.  
One method to leverage both temporal correlations and sparsity is to take a batch processing approach (e.g. ~\cite{GIA:2009,GIA:2011, ZHA:2011,sankaranarayanan2010compressive,ASI:2012mri,aravkin2014optimization,aravkin2013sparse}).
Batch processing is useful in a number of application, in particular when measurements can be amassed over a time period and analyzed simultaneously. 
While batch processing can achieve lower estimation errors and permit greater freedom in cost function design, these methods cannot be run causally, excluding them, for instance, from real-time applications. 
Additionally, operating on larger batches of measurements may require increased computational complexity (e.g. memory or computational time) than online methods.
Thus, batch processing techniques are significantly different from our focus on causal (online) estimation. 
Recent work has also addressed sparse state estimation via efficient updates of previous solutions using specific optimization algorithms, such as homotopy methods~\cite{SAL:2010,SAL:2011} or iterative thresholding~\cite{aureleISTA}. While utilizing the dynamic information implicitly to reduce computational complexity, these approaches do not attempt to improve the estimator quality through modeling signal dynamics. 

Another branch of previous work essentially modifies the Kalman filter equations to induce sparsity in the estimate.  One approach in this vein employs a pseudo-norm in the Kalman update equations to encourage sparser solutions~\cite{KAN:2010} and then enforces an \lo constraint on the state after the Kalman update.   Another approach~\cite{VAS:2008,VAS:2011} works in two-steps, first performing support estimation using an \lo cost function and then running the traditional Kalman equations on a restricted support set.   While the work in~\cite{VAS:2008,VAS:2011} assumes sparse innovations, the robustness to model mismatch has not been fully explored. As will be shown in Section~\ref{sec:results}, these approaches appear to lack robustness to the innovations statistics, perhaps due to the support set estimation being particularly sensitive to estimation errors.  From a computational perspective, storing and inverting full covariance matrices is also prohibitive for the high-dimensional data where sparsity models have been most successful.  

More recent approaches have considered direct coefficient transition modeling via Markov models~\cite{ZAC:2012pursuit,ZIN:2010}, either by utilizing message passing to propagate support information through time~\cite{ZIN:2010} or by using the previous estimate to influence coefficient selection through modified orthogonal matching pursuit~\cite{ZAC:2012pursuit}. These approaches do not rely on solving complex optimization problems and so are extremely computationally efficient.  However, these approaches either restrict the dynamics model to specific applications~\cite{ZIN:2010} (i.e. the approach as specified and implemented restricts the dynamics function to $f(\bm{x}) = \bm{x}$) or have limited robustness due to the fact that they strictly enforce a support set estimate~\cite{ZAC:2012pursuit}.  Additionally, they currently lack strong accuracy guarantees. 

In another direction, some proposed algorithms (e.g.~\cite{VAS:2010,VAS:2010mod,sejdinovic2010bayesian,CHA:2011b,hall2013dynamical}) essentially modify the BPDN objective function to include a term for the dynamic state prediction, including the BPDN-DF algorithm analyzed in the next section.  Such approaches are simple and can work well, but lack robustness to model mismatch due to their explicit strong assumptions on the innovations statistics (see the discussion in~\cite{CHA:2011b}).  Until now, these approaches have also lacked strong accuracy guarantees.




\section{Basis Pursuit De-Noising with Dynamic Filtering}

\label{ssec:BPDNres}

The first algorithm we explore is \emph{basis-pursuit de-noising dynamic filtering} (BPDN-DF)~\cite{CHA:2011b,sejdinovic2010bayesian} (similar to~\cite{VAS:2010,GIA:2011}), which is a simple and efficient method that can work well when the dynamics are known.  BPDN-DF modifies the BPDN estimator to include a term from the Kalman filter optimization that enforces consistency with the model prediction from the previous state estimate
\begin{gather}
	\coefvest_n = \arg\min_{\coefvec} \norm{\insigt{n} - \fvecmat\basemat\coefvec}_2^2 + \gamma\norm{\coefvec}_1 \nonumber \\
	\qquad\qquad + \kthresh\norm{\bm{W}\coefvec - f\left(\statevest_{n-1}\right)}_2^2, \label{eqn:BPDN_DF}
\end{gather}
where $\gamma$ and \kthresh are tradeoff parameters and the state is recovered via $\statevest_n = \basemat\coefvest_n$.  
Thus BPDN-DF seeks a sparse vector that matches both the measurements and the causal prediction obtained via the known dynamics model. 
As depicted in Figure~\ref{fig:Style1}(c), BPDN-DF combines the penalty for deviations from the dynamic prediction $ f\left(\statevest_{n-1}\right)$ with the typical measurement likelihood and the sparsity penalties from BPDN to obtain an overall cost function for estimating the current state.

BPDN-DF has shown empirical performance improvements over standard Kalman filters and independent BPDN estimation at each time-step (i.e., using BPDN at each time step without temporal regularization)~\cite{CHA:2011b}.  Additionally, the algorithm is simple and requires a similar computational complexity as standard BPDN, making it appropriate for high-dimensional data.  However, while BPDN-DF does pass first-order information forward from the previous state to the current estimator, this procedure creates a combined regularizer with a ``rounded'' shape in some directions (illustrated in Figure~\ref{fig:Style1}(c)).  This bulging of the combined prior starts to resemble the $\ell_2$ ball, and is in conflict with the geometric structures known to produce sparse solutions (e.g., the $\ell_1$ ball)~\cite{BAR:2010PIEEE,BAR:2007}.  This can cause robustness problems in BPDN-DF when there are inaccuracies in the dynamics model (demonstrated later).  

Aware of these potential drawbacks, analytic performance guarantees would be valuable because BPDN-DF remains an effective and computationally attractive algorithm for some applications. Prior work in~\cite{hall2013dynamical} provided guarantees on the statistical regret of an estimator similar to BPDN-DF called the dynamic mirror descent (DMD). Although the DMD is more general in terms of its formulation, the guarantees it provides (i.e., bounded regret) are in terms of a comparison with an optimal non-dynamical cost function and are therefore difficult to directly relate to the cost function in~Equation~\ref{eqn:BPDN_DF}. While our preliminary analysis of BPDN-DF using simple techniques provided some accuracy guarantees, these bounds were too loose to be an accurate reflection of the estimator performance~\cite{charles2014globalsip1}. In this work we adopt a proof technique similar to~\cite{aureleISTA,bredies2008linear} which provides much stronger accuracy guarantees for BPDN-DF, as summarized in the following theorem.

\begin{thm}
    \label{thm:BPDNDFbasic}
    Suppose that at each time $n$, $\mematt{n}$ satisfies RIP($S+2q$,\ripd) for some constant $q\geq0$, and the error and innovations satisfy $\norm{\meerrvt{n}}_2\leq\errmax$ and $\norm{\innovvt{n}}_2\leq\inmax$. Furthermore, suppose that for all $\bm{x}_1,\bm{x}_2$, the dynamics satisfies $\norm{f\left(\bm{x}_1\right) - f\left(\bm{x}_2\right)}_2 \leq f^{\ast}\norm{\bm{x}_1 - \bm{x}_2}_2$ for a universal constant $f^{\ast}$, and the coefficients at each $n$ satisfy $\|\bm{z}_n\|_2 \leq b$.  If $\gamma > 0$, $\kthresh>0$ are known constants and the following condition is met:
    \begin{gather}
        \kthresh\left((1 - \dlip)b - (1+\ripd)\gamma\sqrt{q} - \inmax \right) \geq \qquad\qquad \nonumber \\
	\qquad\qquad (1 + \ripd)\gamma\sqrt{q} + \sqrt{1+\ripd}\errmax - \left(1 - \ripd \right)b, \nonumber
    \end{gather}
    then the solution to Equation~\eqref{eqn:BPDN_DF} satisfies  
    \begin{gather}
	    \norm{\coefvtest{n} - \coefvect{n}}_2 \leq  \beta^n\norm{\errvt{0}}_2 + \left(1-\beta^n\right)\left(C_1\sqrt{q} + C_2\errmax + C_3\inmax  \right) \nonumber
    \end{gather}
    where the constants are defined as:
    \begin{eqnarray}
        \beta & = & \frac{\kthresh\dlip}{1+\kthresh - \ripd} \nonumber \\
        C_1 & = & \frac{(1+\ripd)(1+\kthresh)\gamma}{1 + \kthresh - \ripd - \kthresh\dlip} \nonumber \\
        C_2 & = & \frac{\sqrt{1+\ripd}}{1 + \kthresh - \ripd - \kthresh\dlip} \nonumber \\
        C_3 & = & \frac{\kthresh}{1 + \kthresh - \ripd - \kthresh\dlip}.
    \end{eqnarray}
\end{thm}

The full proof of Theorem~\ref{thm:BPDNDFbasic} is outlined in Appendix~\ref{app:bpdndfproof}. 
Note that the error bound has two distinct components: a transient term depending on the initial uncertainty that decays if $\beta^n\rightarrow 0$, and a steady state term that emerges depending on the inherent difficulty of the problem (i.e., the noise sources and the sparsity of the signal).  If $\kthresh = 0$ (removing the dynamics in the estimation), we obtain exactly the results in~\cite{aureleISTA} for solving BPDN with no dynamic filtering term. 
Alternatively, as $\gamma \rightarrow 0$, the result in Theorem~\ref{thm:BPDNDFbasic} ceases to hold as convergence is no longer guaranteed. Specifically, $\gamma = 0$ violates a condition in Lemma~\ref{lem:ubound}, which is required to prove algorithmic convergence. 
We note that the same condition in Lemma~\ref{lem:ubound} necessitates the condition $\|\bm{z}_n\|_2 \leq b$, indicating that $b\rightarrow \infty$ can also cause a loss of convergence, and thus accuracy, guarantees. 

Theorem~\ref{thm:BPDNDFbasic} implies that BPDN-DF is only guaranteed to converge if $\beta < 1$. Since most parameters are system or signal dependent (and therefore not controllable by the user), we interpret this as a condition on the parameters $\gamma$ and \kthresh. In particular, since $\gamma$ does not appear in the expression for $\beta$, we can consider this to be a bound on \kthresh,
\begin{gather}
	\kthresh  <  \frac{1 - \ripd}{\dlip - 1}  \quad \mbox{if} \quad \dlip > 1.    \label{eqn:kcond}
\end{gather}
This condition essentially compares the smoothness of the dynamics function with the RIP conditioning of the measurements. For example, as $\ripd$ approaches 1, the allowable range of \kthresh shrinks indicating that the dynamics should be deemphasized in BPDN-DF. Alternatively, as $\dlip$ approaches $1$,  the range of \kthresh increases, indicating that the dynamics can be emphasized in the optimization cost. Interestingly, if $\dlip \leq 1$, the requirement that $\beta < 1$ induces no additional restrictions on \kthresh as the numerator of Equation~\eqref{eqn:kcond} is positive.  With a negative denominator in Equation~\eqref{eqn:kcond}, this condition is redundant with the prior requirement that \kthresh is positive.  
The steady-state error bound depends on both the maximum measurement error energy \errmax and the maximum innovations energy \inmax. The  parameter \kthresh trades off between the two, where the trade-off takes into account the RIP conditioning \ripd as well as the dynamics smoothness \dlip.





\section{Re-Weighted \lo Dynamic Filtering}

\label{ssec:ModUp}

\begin{figure}[t]
	\centering
	\includegraphics[width=2in]{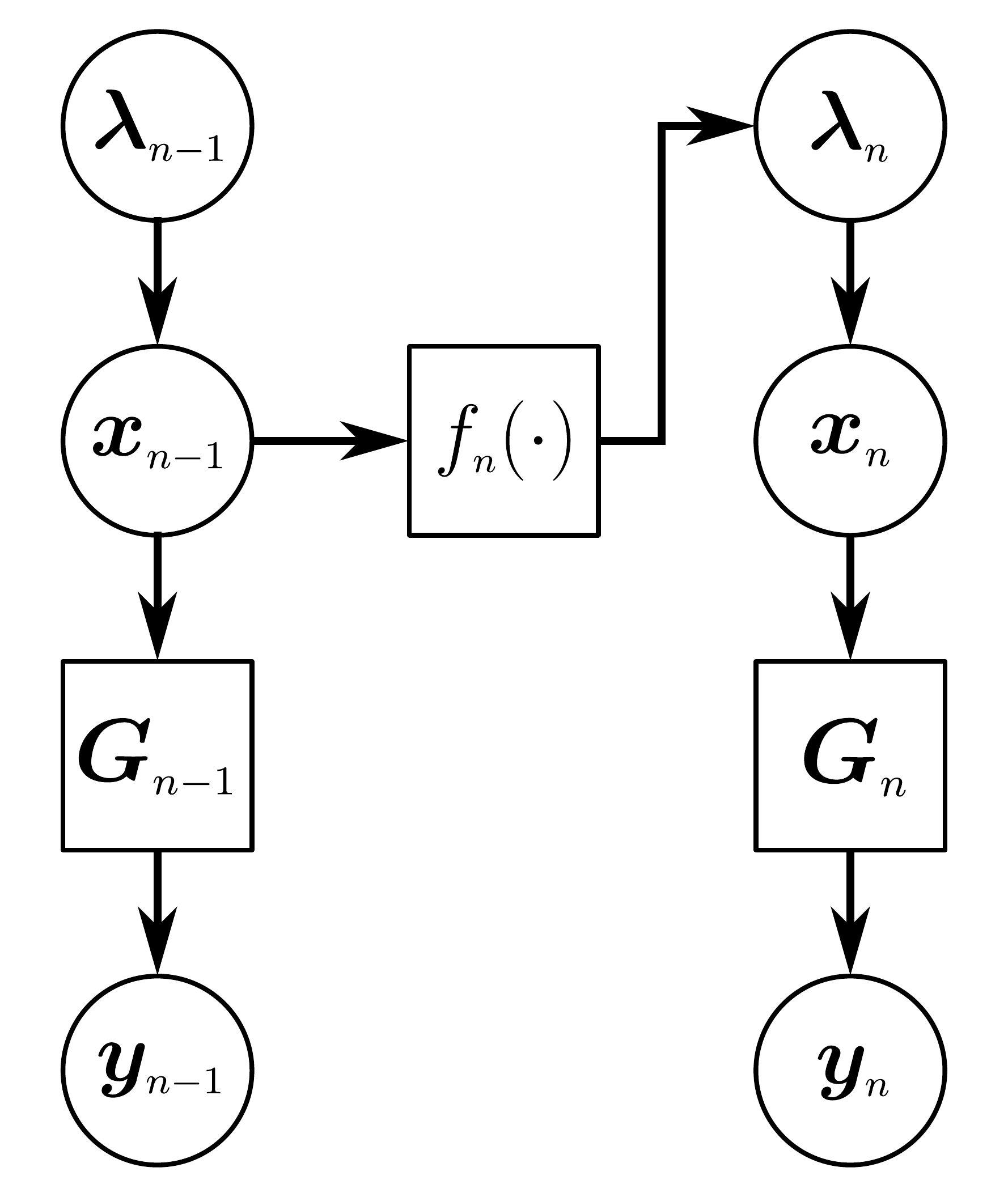}
	\caption{The RWL1-DF algorithm inserts the dynamic information at the second layer of the LSM model for each time step. The graphical model depicts the model dependencies, where prior state estimates are used to set the hyperpriors for the second level variables controlling the variances (i.e., SNRs) of the state estimates at the next time step.} 
	\label{fig:rwl1_gmod}
\end{figure}

Despite now having provable guarantees, performance of BPDN-DF degrades with increased inaccuracy in the dynamics model. 
In fact, Theorem~\ref{thm:BPDNDFbasic} requires that the innovations (a measure of the dynamics inaccuracy) satisfies certain inequalities with respect to the measurement noise and sparsity. 
In this section we present a novel algorithm called \emph{re-weighted $\lo$ dynamic filtering} (RWL1-DF) that uses a richer statistical model to improve performance, especially when the dynamics model is inaccurate.  
In contrast to existing algorithms (including BPDN-DF), RWL1-DF displays a combination of flexibility (allowing arbitrary dynamic models), computational efficiency (leveraging advances in \lo minimization), and robustness to dynamic model errors caused by sparse innovations. 
The main technical innovation we propose is to use the hyper-priors in the LSM model to incorporate a dynamic prediction into second-order statistics for the next estimate in a causal fashion (akin to Kalman filtering).  
The resulting estimation procedure is illustrated in Figure~\ref{fig:Style1}(d) and the graphical model is depicted in Figure~\ref{fig:rwl1_gmod}.   While some existing literature~\cite{HAS:2009,VAS:2010mod,mio2009comp,car2010it} has proposed re-weighted estimation schemes to include prior information in a limited way (e.g., using only binary weights based on support set estimates), RWL1-DF modifies RWL1 specifically for dynamic filtering and uses a general probabilistic model that improves robustness by allowing arbitrary weightings.

Technically, similar to the LSM model, our proposed model describes the conditional distribution on the sparse coefficients $\bm{z}$ at time $n$ as a zero-mean Laplacian with different variances: 
\begin{gather}
	p\left(\bm{z}_n[i] \left|\right. \threshvec_n\right)  = \frac{\lambda_0\threshvec_n[i]}{2}e^{-\lambda_0\threshvec_n[i]\left|\bm{z}_n[i]\right|}. \nonumber
\end{gather}
We also set the scale variables controlling the coefficient variance to be Gamma distributed: 
\begin{gather}
	p\left(\threshvec_n[i]\right) =  \frac{\threshvec_n^{\alpha-1}[i]}{\bm{\theta}^{\alpha}_n[i]\Gamma\left(\alpha\right)}e^{-\threshvec_n[i]/\bm{\theta}_n[i]}, \nonumber
\end{gather}
and allow each of these variables to have a different expected value, ($\alpha\bm{\theta}_n[i] = E[\threshvec_n[i]]$) by modifying  $\bm{\theta}_n[i]$.  
To introduce dynamic information, the expected value of each scale variable is set based on the previous state estimate.  In particular, if the model prediction of $\bm{z}_n[i]$ based on the previous state is large (or small), the variance of that coefficient at the current estimate is made large (or small) by making $\threshvec_n[i]$ small (or large).  
Large variances allow the model flexibility to choose from a wide-range of non-zero values for coefficients that are likely to be active and small variances (with a mean of zero) encourage the model to make such coefficients inactive.  By ``encouraging'' the model through the use of second order statistics (instead of forcing the model to use an estimated coefficient subset) the model remains robust and flexible.
While there are many possible choices for the mapping of predictions into the model at the next time step, in this paper we choose 
\[ \bm{\theta}_n[i] = \xi\left(|\bm{W}^{-1}\dfunct{\bm{W}\widehat{\bm{z}}_{n-1}}{n}[i]| + \eta\right)^{-1}, \]
where $\bm{W}^{-1}\dfunct{\bm{W}\widehat{\bm{z}}_{n-1}}{n}[i]$ is the $i^{\mbox{th}}$ coefficient of the previous signal propagated through the dynamics, $\eta$ is a linear offset and $\xi$ is a multiplicative constant. Note that any general model for the dynamics is allowable and we are not restricted to linear dynamical systems in the prediction.  Note also that the absolute values are necessary because \threshvec determines a variance, which must be strictly positive.  The parameter $\eta$ determines the distribution of the variance when the coefficient is predicted to be zero, resulting in $\bm{W}^{-1}\dfunct{\bm{W}\widehat{\bm{z}}_{n-1}}{n}[i] = 0$.  This parameter reflects the magnitude of the innovations in the erroneous model predictions.

Under this model, the joint MAP estimate of all model parameters becomes
\begin{eqnarray}
	[\widehat{\bm{z}}_n, \widehat{\threshvec}_n] & = & \arg\min_{[\bm{z}, \threshvec]} \norm{\insigt{n} - \mematt{n}\bm{W}\bm{z}}_2^2 + \beta\sum_i{\left|\threshvec[i]\bm{z}[i]\right|} \nonumber \\
	& & - \alpha\sum_i{\log\left(\threshvec[i]\right)} \label{eqn:RWL1DF_est} \\
	& & + \frac{\lambda_0}{\xi}\sum_i{\threshvec[i]\left(|\bm{W}^{-1}\dfunct{\bm{W}\widehat{\bm{z}}_{n-1}}{n}[i]| + \eta\right)}. \nonumber
\end{eqnarray} 
As in the LSM case, the optimization in~\eqref{eqn:RWL1DF_est} is not easily solved for both $\bm{z}_n$ and $\threshvec_n$ jointly, but the model yields a simple form when using an EM approach.
The precise steps of the iteration in this case are 
\begin{eqnarray}
	\mbox{E step:}	& \quad	&	\widehat{\threshvec}^t_n = E_{p\left(\threshvec {|} \widehat{\bm{z}}_n^t\right)}\left[\threshvec \right]  \nonumber \\
	\mbox{M step:}	& \quad	&	\widehat{\bm{z}}_n^t = \arg\min_{\bm{z}_n}-\log\left[p\left(\bm{z}_n{|}\widehat{\threshvec}\right)\right] \nonumber
\end{eqnarray}
where $t$ denotes the EM iteration number and $E_{p\left(\threshvec {|} \widehat{\bm{z}}_n\right)}\left[\cdot\right]$ denotes the expectation with respect to the conditional distribution $p\left(\threshvec {|} \widehat{\bm{z}}_n\right)$. We can write the maximization step as 
\begin{gather}
	\widehat{\bm{z}}_n^t = \arg\min_{\bm{z}} \norm{\insigt{n} - \memat\bm{W}\bm{z}}_2^2 +  \lambda_0\sum_i{\left|\widehat{\threshvec}^t[i]\bm{z}[i]\right|}, \label{eqn:drwl1_opt1}
\end{gather}
since the MAP optimization conditioned on the \thresh parameters reduces to a weighted \lo optimization. 

The conjugacy of the Gamma and Laplacian distributions admit a simple closed-form solution for the expectation step in this model that is similar to the original LSM model~\cite{GAR:2010},
\begin{gather}
	\threshvec_n^{t+1}[i] = \frac{2\tau}{\beta|\widehat{\bm{z}}^t_n[i]| +  |\bm{W}^{-1}\dfunct{\bm{W}\widehat{\bm{z}}_{n-1}}{n}[i]| + \eta}, \label{eqn:drwl1_opt2}
\end{gather}
where $\tau = \left(\alpha+1\right)\xi$ is a constant scaling value, $\beta = \lambda_0\xi$ can be interpreted as a tradeoff between the measurement and the prediction, and the signal of interest can again be recovered via $\statevestt{n} = \bm{W}\widehat{\bm{z}}_n$. 
The resulting procedure that iterates between Equations~\eqref{eqn:drwl1_opt1} and~\eqref{eqn:drwl1_opt2} looks nearly identical to the static RWL1 algorithm except that the denominator in the \threshvec update contains a term depending on the previous state. This term encourages smaller \threshvec values (i.e., higher variances) in the elements that are predicted to be highly active according to \dfunct{\statevtest{n-1}}{n}. This graduated encouragement of coefficients selected by the prediction allows the algorithm to perform especially well when the states and innovations are sparse while retaining good performance when the innovations are denser. Furthermore, the simple form means that the recent work in \lo optimization methods can be leveraged for computationally efficient recursive updates, requiring only a small number of \lo optimizations at each time step.  In particular, since no covariance matrix inversion is required and many modern \lo estimation methods require only matrix multiplication (and no inversion), this approach is also amenable to high-dimensional data analysis.

Despite being highly nonlinear, RWL1-DF has demonstrable stability and convergence properties.  First, the RWL1-DF update equations in~\eqref{eqn:drwl1_opt1} and~\eqref{eqn:drwl1_opt2} directly guarantee that the coefficient and variance estimates will be bounded.  To see this for an arbitrary time step $n$, note that the variance estimates are within the range $\threshvec_n^{t+1}[i]\in(0, 2\tau/\eta]$ for each EM iteration. With bounded variances, the weighted BPDN optimization in~\eqref{eqn:drwl1_opt2} will yield a solution where the output coefficients are also bounded. 
Second, existing guarantees for the EM algorithm guarantee that the EM iterations in the proposed algorithm (i.e., an estimate at a single time step) have coefficient differences that asymptotically converge to zero, $\lim_{t\rightarrow\infty}\|\bm{z}_n^t -\bm{z}_n^{t+1}\| = 0$~\cite{chret2000kullback}.  Existing stronger convergence guarantees (from EM or from the very limited analysis on RWL1 algorithms~\cite{needell2009noisy}) do not apply in this case because of the technical details of RWL1-DF.
Despite the lack of stronger convergence guarantees, the numerical results in section~\ref{sec:results} demonstrate that the algorithm converges in practice with just a few EM iterations.




\section{Simulation Results}
\label{sec:results}

While the dynamic filtering approaches discussed here are general inference tools with many potential applications, we focus our  evaluation on compressed sensing (CS) recovery as an example test case.  In all simulations we compare the performance with existing algorithms discussed in Section~\ref{sec:SKFback} when possible, noting in each particular case where algorithms were unable to be evaluated because they are incompatible or computationally prohibitive.  
Standard Kalman filtering is not shown because it performs very poorly in these type of simulations (i.e., it doesn't converge to a stable estimate with the sparse statistics of the applications we use)~\cite{CHA:2011b}.   
The performance of independent BPDN (Equation~\eqref{eq:l1optreg}) and independent RWL1 (Equations~\eqref{eqn:RWL1reg}), each applied independently at each time step with no temporal information, are also shown to highlight the benefit of including dynamic information.

\subsection{Stylized tracking scenario}

To explore the performance and robustness of the algorithms in detail, we first perform inference on synthetic data that simulates a stylized tracking scenario.  The use of synthetic data provides us with ``ground truth'' so we can make controlled variations of the data characteristics.  
In this setup, we generate an image with $S$ non-zero moving pixels of various intensities that move with time and represent targets that must be tracked.  The movement of these non-zero pixels $\bm{F}_k$ is specified to be constant motion, and the simulated dynamics includes  a sparse innovations term (i.e., dynamic model error) that causes target motion to change in each time step for some percentage of the pixels $p$.  In other words, at every time step there is a probability $p$ of each target abruptly changing directions to violate the dynamics model assumed by the inference algorithms.  This process simulates an innovation that is approximately $2Sp$-sparse at every iteration,  allowing us to evaluate the algorithm's robustness to a type of model mismatch (i.e., shot noise) that is particularly challenging for Kalman filter techniques.

In detail, we create 24x24 pixel videos ($N = 24^2$) with 20 moving particles ($S = 20$). The vectors are observed with Gaussian measurement matrices (with normalized columns) that are independently drawn at each iteration, and we add Gaussian measurement noise with variance $\sigma_{\bm{\epsilon}}^2 = 0.001$. We vary the number of measurements to observe the reconstruction capability of the algorithm in highly undersampled regimes, but the number of measurements per time step is always constant within a trial.   All simulations average the results of 40 independent runs and display reconstruction results as the relative mean-squared error (rMSE) for each frame, calculated as:
\begin{gather}
	\frac{\norm{\statev_k - \statevest_k}_2^2}{\norm{\statev}_2^2}.
\end{gather}
For independent BPDN, at each iteration we use the value $\lambda$ = 0.55$\sigma_{\bm{\epsilon}}^2$. For independent RWL1 we use $\lambda_0$ = 0.0011, $\tau$ = 1 and $\nu$ = 0.01.  For BPDN-DF we use $\gamma$ = 0.5$\sigma_{\bm{\epsilon}}^2$ and $\kappa$ = $0.0007/(p+1)$. For RWL1-DF we use $\lambda_0$ = 0.0011, $\tau$ = 1 and $\nu$ = $1 - 2p/\sparsity$. These parameters were optimized using a manual parameter sweep.   Furthermore, for comparison we also show the performance of an optimal oracle least-squares solution, where the support at each iteration is known ($\mematt{k}$ becomes overdetermined).  Note that this scenario is particularly challenging for many existing algorithms because of the arbitrary model dynamics (e.g., $\bm{F}_k\neq I$) that may vary with time.  In particular, comparisons with the algorithms described in~\cite{VAS:2010mod,VAS:2010,mio2009comp,HAS:2009} (or modifications of them to accommodate arbitrary dynamics) were attempted and are not shown here because they still performed significantly worse than static estimation (e.g., BPDN) even after extensive searching for good parameter settings.

\begin{figure}[t]
	\begin{center}
		\includegraphics[width=0.4\textwidth]{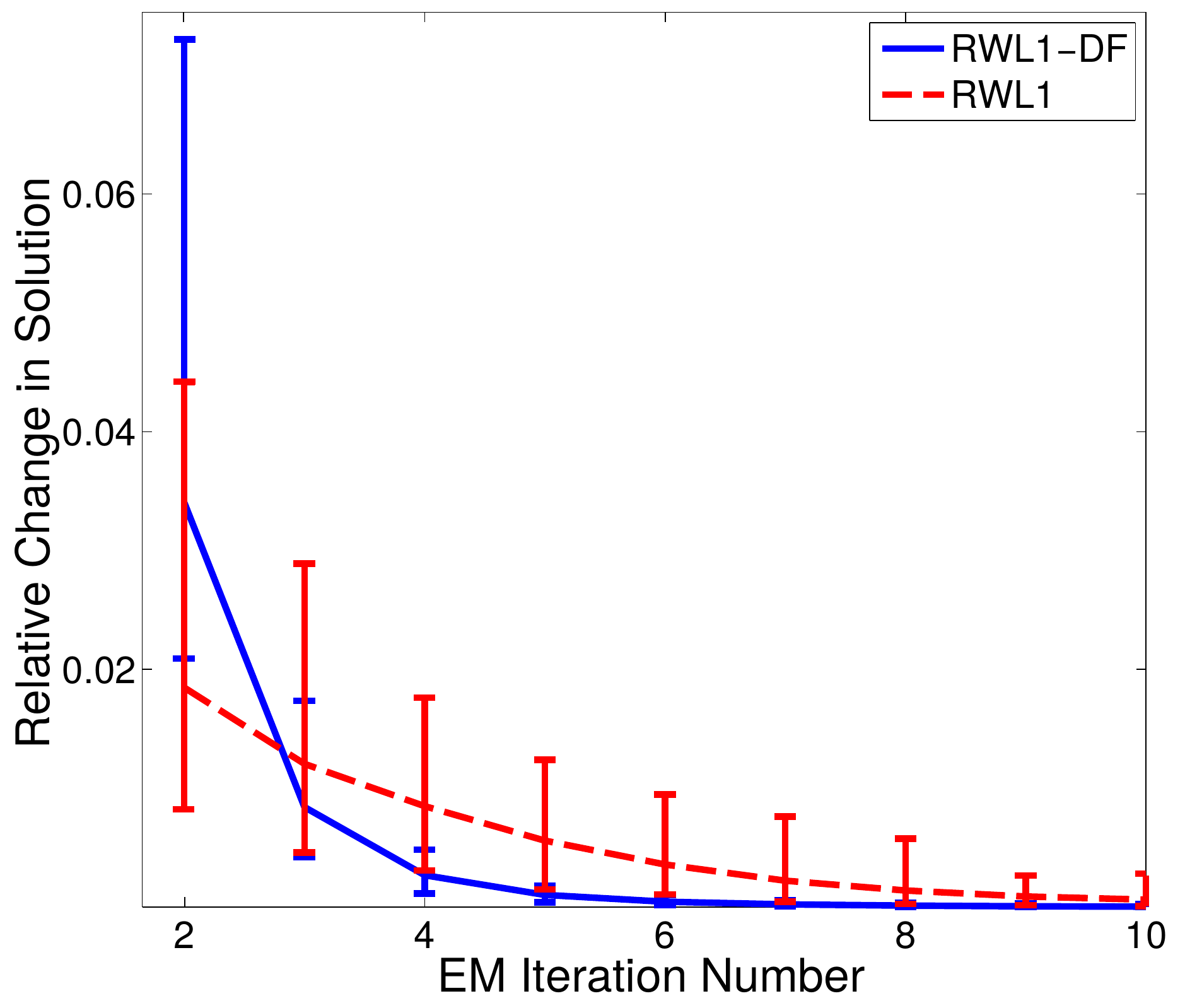}
	\end{center}
	\caption{Convergence of RWL1 and RWL1-DF. The mean relative change (over 200 frames) in the coefficients is plotted, with the error bars indicating the maximum and minimum values. The relative norm difference of the coefficients in the RWL1 and RWL1-DF algorithms falls quickly over the first 10 iterations. The dynamic information helps the RWL1-DF algorithm converge faster, requiring approximately 5-7 iterations to converge. }
	\label{fig:RWL1conv}
\end{figure}

\begin{figure*}[t]
	\centering
	\includegraphics[width=0.9\textwidth]{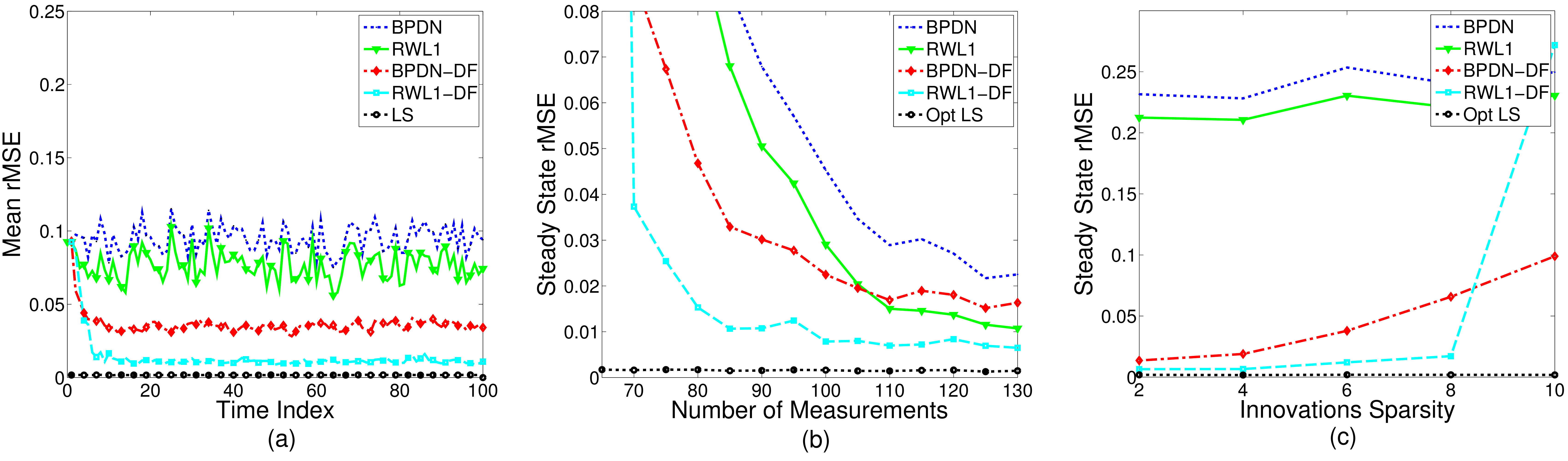}
	\caption{Behavior of the RWL1-DF algorithm on synthetic data. (a) RWL1-DF converges  to a lower mean rMSE than static sparse estimation or BPDN-DF.  Shown for \measdim = 80, \statedim = 576, $S$ = 20 and $p$ = 0.25. (b) When sweeping the number of measurements \measdim for \statedim = 576, $S$ = 20, and $p$ = 0.25, we  observe that the performance improvement for RWL1-DF is especially distinct in the highly undersampled regime.  Each point is the average steady state rMSE over 40 independent trials. (c) RWL1-DF is also more robust to model mismatch in the innovation statistics.  Shown here for different innovations sparsity ($2Sp$) for \measdim = 70, \statedim = 576, and $S$ = 20. Each data point is the result of averaging the steady-state rMSE over 40 independent trials. Note that when BPDN-DF starts to perform better ($2Sp=10$), the innovations are actually half of the total support set and a Gaussian innovations model may be more accurate than a sparse innovation model. 
	}
	\label{fig:tsweep}
\end{figure*}

The iterative algorithms based on RWL1 are stopped when the relative norm-squared difference between coefficients at consecutive iterations $\|\bm{z}_n^t - \bm{z}_n^{t+1}\|_2^2/\|\bm{z}_n^t\|_2^2$ falls below a specified threshold (we use 0.1\%), which typically requires just a few (typically 5-7) iterations. Figure~\ref{fig:RWL1conv} shows the relative coefficient change over EM iteration, demonstrating that RWL1 (i.e. without dynamics) convergence occurs by 10 iterations and RWL1-DF actually converges faster due to the improved performance from incorporating dynamic information.
Figure~\ref{fig:tsweep}(a) shows a single trial with \measdim = 80 measurements and $2Sp$ = 5 innovation errors at each time step.  Note that RWL1-DF reaches substantially lower steady-state recovery error than BPDN-DF, illustrating the net improvements gained by using second-order statistics in the estimation.  Figure~\ref{fig:tsweep}(b) displays the results of varying the number of measurements while holding $2Sp = 5$. While performance for all algorithms becomes comparable for large numbers of measurements, it is clear that exploiting temporal information can most improve performance in the highly undersampled regime.  In particular, RWL1-DF is able to sustain virtually the same steady-state rMSE down to much more aggressive levels of undersampling than BPDN-DF.  Finally, we explore the robustness of each algorithm to model errors by fixing the number of measurements ($\measdim = 70$) and varying the sparsity of the innovations $2Sp$.  Figure~\ref{fig:tsweep}(c) shows the results, illustrating that RWL1-DF uses the second-order statistics to sustain better performance than BPDN-DF when the innovations are sparse (i.e., shot noise).  
We note that when $2Sp>8$, the total number of model errors is 50\% of the signal sparsity, indicating both that the innovation is dense, and that the innovations has as much, or more, energy than the signal. The combination of the innovations statistics mismatch (sparse vs. dense) and reduced predictive power of the dynamics model produces a phase transition where RWL1-DF is no longer able to correct large innovation values.

\subsection{CS recovery of natural video sequences}

To test the utility of RWL1-DF on natural signals, we explore its performance on a simulation of compressively sampled natural video sequences.  These results will report in-depth comparison of a single challenging video sequence (the Foreman sequence\footnote{The Foreman sequence is available at: http://www.hlevkin.com/TestVideo/foreman.yuv}) as well as aggregate statistics from a batch of video from a BBC nature documentary (as used in~\cite{cad2012learning}).  The documentary footage is valuable as broad comparison because it contains many different types of motion, including static frames with localized changes and highly dynamic frames with moving subjects across large portions of the visual field.

In our simulation of CS video recovery, we take the time-varying hidden state $\statevt{n}$ to be the wavelet (synthesis) coefficients at each frame of the video.
While the true frame-by-frame dynamics of natural video are likely to be complex and non-linear, for this simulation we use a simple first-order model predicting that the coefficients will remain the same from one frame to the next: \dfunct{\statev}{k} = \statev for all $k$.  Improvements to this simple model would improve any algorithm that can incorporate more advanced models, including RWL1-DF which allows arbitrary dynamics models.  A benefit to assuming stationary dynamics is that it allows us to compare recovery performance with a number of existing algorithms that do not currently have arbitrary dynamics as part of the approach, including DCS-AMP~\cite{ziniel2012dynamic} and modCS~\cite{VAS:2010mod}. To demonstrate the advantage of using the coefficient values instead of only support information, we also compare to a modified version of the binary-weighted BPDN algorithm described in~\cite{HAS:2009} where we have added support set prediction based on model dynamics.  We call this approach weighted $\ell_1$ with prior information (WL1P).

\begin{figure*}[t]
	\centering 
	\includegraphics[width=\textwidth]{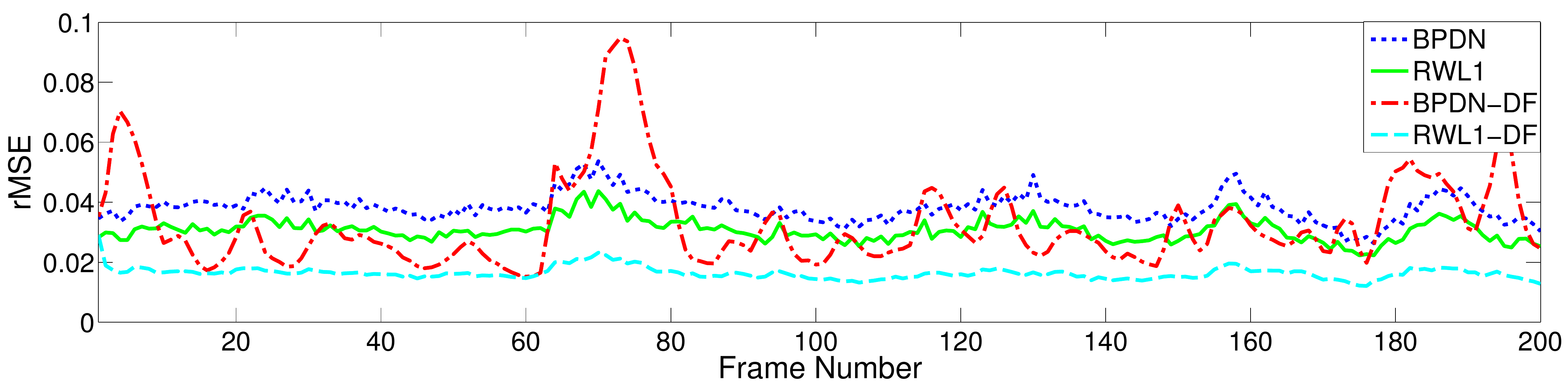} 
        \caption{CS recovery of the full Foreman video sequence. Each curve represents the rMSE for recovery from subsampled noiselets ($M$/$N$ = 0.25) using the DT-DWT as the sparsifying basis. The independent BPDN recovery (dotted blue curve) and the independent re-weighted BPDN (solid green curve) retain a steady rMSE over time. RWL1-DF (dashed cyan curve) converges on a lower rMSE than either time-independent estimation and remains at approximately steady-state for the remainder of the video sequence. BPDN-DF (the dot-dash red curve) can converge to low rMSE values, but is highly unstable and can yield very poor results when the model is not accurate due to motion in the scene.} 
	\label{fig:vidtrial}
\end{figure*}

We use a four-times overcomplete  dual-tree discrete wavelet transform (DT-DWT)~\cite{selesnick2005dual} as the sparsity basis and we simulate CS measurements by applying a subsampled noiselet transform~\cite{coif2001noiselet} to each frame.  We take $\measdim$ = 0.25\statedim measurements per frame (where \statedim = $128^2$) with a measurement noise variance of $10^{-4}$.\footnote{We use the noiselet transform because it can be computed with an efficient implicit transform and has enough similarity to a random measurement that it works well in CS~\cite{romberg2008imaging}.}  We solve all optimization programs using the TFOCS package~\cite{BEC:2011} due its stability during RWL1 optimization and the ability to use fast implicit operators for matrix multiplications. 
In the Foreman video sequence, we simulate CS measurements on a 128$\times$128 pixel portion of the video, using parameters  
$\lambda$ = 0.001 for BPDN, $\lambda_0$ = 0.11, $\tau$ = 0.2 and $\nu$ = 0.1 for RWL1, $\gamma$ = 0.01 and $\kappa$ = 0.2 for BPDN-DF, and $\lambda_0$ = 0.003, $\tau$ = 4, $\beta$ = 1 and $\nu$ = 1 for RWL1-DF. 
These parameters were found using a manual parameter sweep to optimize performance for each algorithm and number of CS measurements.  While the number of measurements was fixed in this simulation for computational tractability, recovery performance could be altered for all algorithms by adjusting the number of measurements.

\begin{figure*}[t]
	\centering
	\includegraphics[width=0.9\textwidth]{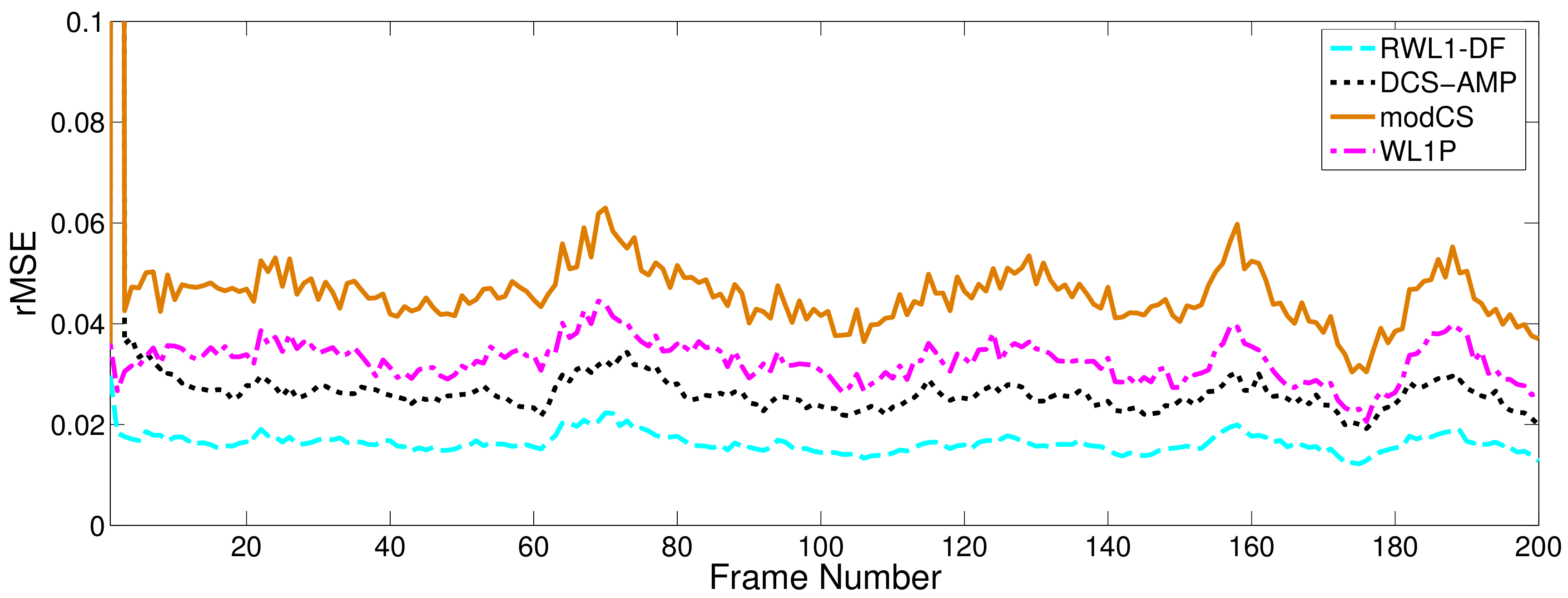}
	\caption{A comparison of RWL1-DF with existing recovery algorithms (DCS-AMP, modCS and WL1P) for the Foreman video sequence. Each curve represents the rMSE for recovery from subsampled noiselets ($M$/$N$ = 0.25) using the DT-DWT as the sparsifying basis.}
	\label{fig:vidtrial2}
\end{figure*}

\begin{figure*}[t]
	\centering
	\includegraphics[width=0.8\textwidth]{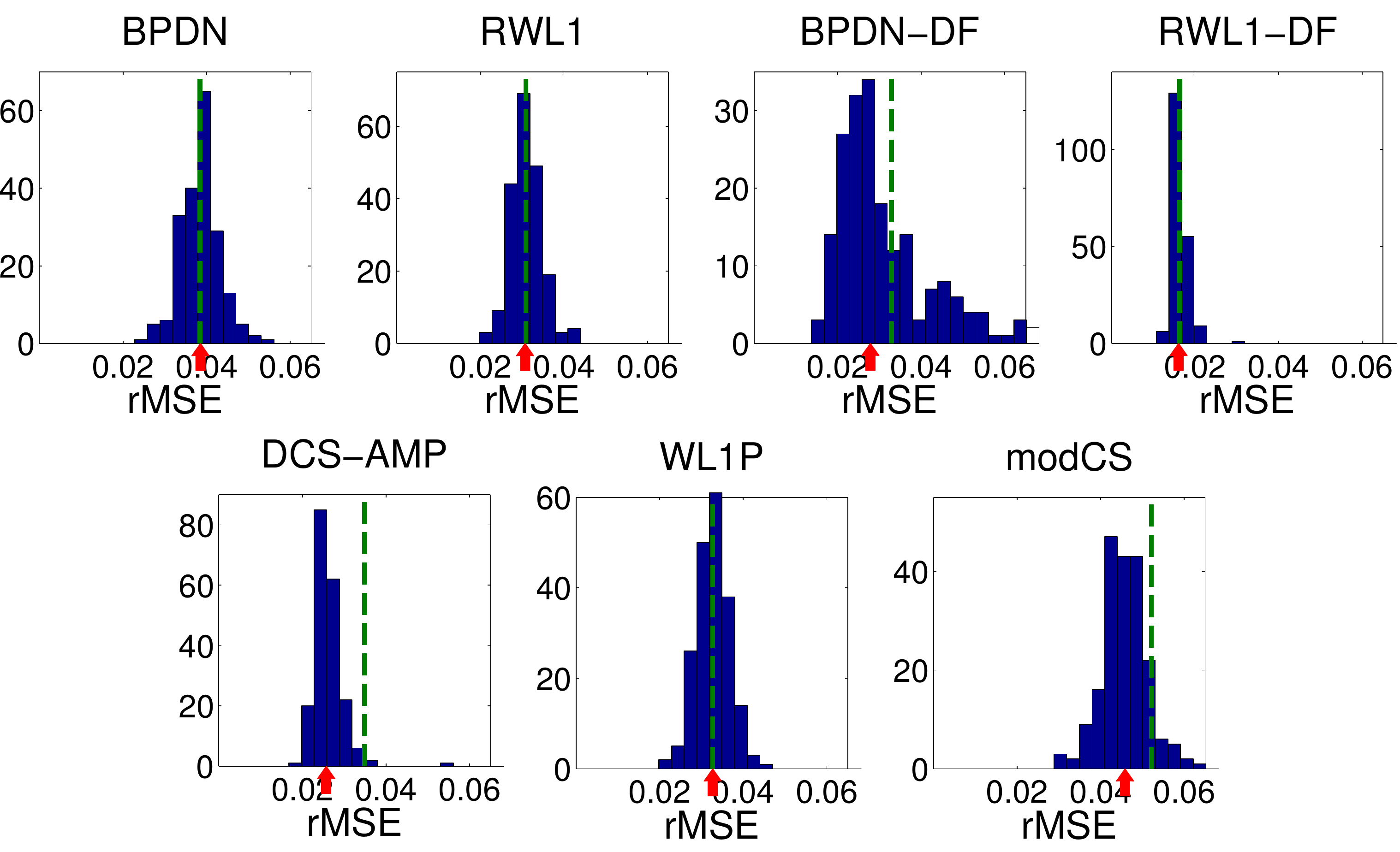}
	\caption{Histogram of the rMSE for the compared algorithms when recovering the Foreman video sequence with the DT-DWT as the sparsifying basis. RWL1-DF  achieves a lower mean (indicated by the dashed green lines) and median (indicated by the red arrows), with a tightly concentrated error distribution due to the robustness to model mismatch (producing few outliers).  Specific mean and median values are shown in Table~\ref{tab:MeanMed}.}
	\label{fig:vidhist}
\end{figure*}

Figure~\ref{fig:vidtrial} shows the recovery of 200 consecutive frames of video in the Foreman sequence. We see that RWL1-DF converges to the lowest steady-state rMSE and is able to largely sustain that performance over the sequence. In contrast, BPDN-DF cycles through periods of good performance and poor performance,  sometimes performing worse than not using temporal information at all.  In essence, BPDN-DF is not robust to model errors, and each time there is motion in the scene (violating the simple dynamics model) the algorithm has to re-converge.  The RWL1-DF approach does not exhibit this performance oscillation because the use of second-order statistics to propagate temporal information is less rigid, allowing for more robustness during model errors.  
Figure~\ref{fig:vidtrial2} shows a comparison of recovery for the Foreman video sequence across several existing algorithms, including DCS-AMP, modCS and WL1P.  Again we optimize algorithms parameters manually to achieve the best aggregate performance over the video sequence. To summarize the performance over the entire video sequence, Figure~\ref{fig:vidhist} shows histogram plots of the rMSE values over the 200 frames for each recovery.  The mean and median for each histogram are represented by the green dashed line and the red arrow respectively, and are listed in Table~\ref{tab:MeanMed}.
The recovery errors for BPDN-DF are significantly spread out, achieving nearly the same median error as the independent algorithms and having a much higher mean error due to the large excursions during model mismatch. BPDN-DF actually may be a worse choice than not using any temporal information at all when model errors are present.  In contrast, RWL1-DF shows a much tighter distribution of errors, having a mean and median significantly lower than alternate approaches.  

\begin{table*}
	\centering
	\begin{tabular}{|l||c|c|c|c|c|c|c|}
		\hline
		              & BPDN    & RWL1  & BPDN-DF  & RWL1-DF &DCS-AMP & WL1P    & modCS    \\ \hline\hline
		Mean rMSE     & 3.84\% & 3.09\% & 3.29\%  &\textbf{1.63\%}  & 3.48\%  & 3.27\%  & 5.22\%   \\ \hline
		Median  rMSE  & 3.85\% & 3.07\% & 2.78\% &\textbf{1.61\%}   & 2.57\%  & 3.27\%  & 4.58\%   \\ \hline
	\end{tabular}
	\caption{Mean and median values for compressive recovery of the Foreman video sequence.}
	\label{tab:MeanMed}
\end{table*}

\begin{figure*}[t]
	\begin{center}
		\includegraphics[width=0.84\textwidth]{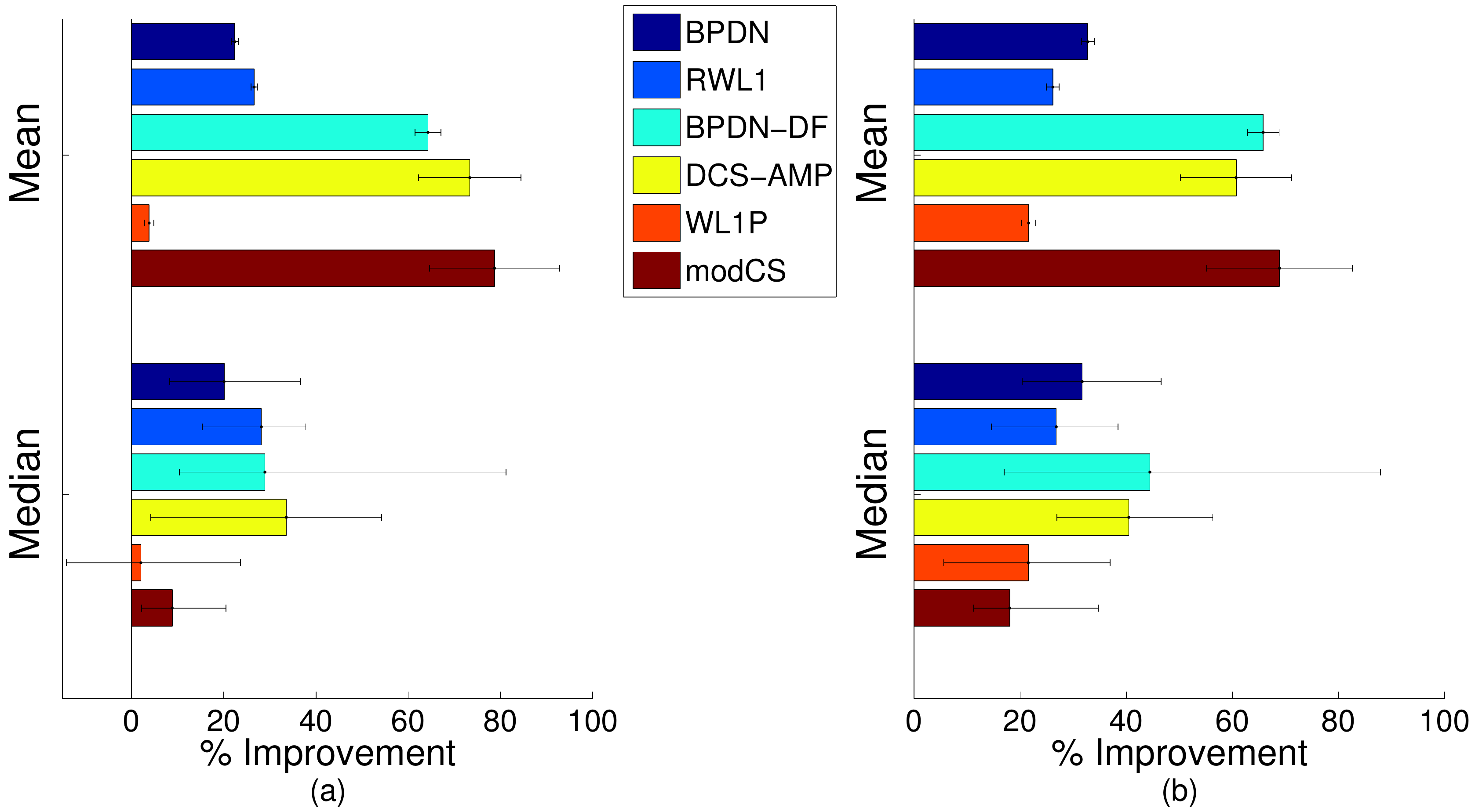}
	\end{center}
	\caption{Percent improvement of RWL1-DF over other algorithms for compressive recovery of video sequences (calculated as described in the text). Displayed is both the mean improvement (error bars indicating the normalized standard deviation) and the median improvement (error bars indicating the 25th and 75th percentile) for each algorithm. (a)  Results for the full database of 24 sequences from a BBC nature documentary. (b) Results for the 13 sequences that were especially challenging for CS recovery, illustrating the benefits of RWL1-DF in this particularly difficult regime.}
	\label{fig:ImpBar}
\end{figure*}

In addition to the in-depth comparison on the single Foreman video sequence, we also perform the same CS recovery task on a database of video sequences from a nature BBC documentary to investigate the performance across a wider range of video characteristics (i.e., including video clips with localized motion and global motion in the scene). We simulated CS measurements for 24 sequences (48-frames each) in the same manner as the Foreman sequence and recovered the frames using the same methodology described above (including parameters optimized on the Foreman sequence). Figure~\ref{fig:ImpBar} shows the mean and median \emph{improvement} of RWL1-DF relative to the other algorithms being evaluated.  Specifically, the plotted mean improvement is the average of the rMSE difference between RWL1-DF and the comparison algorithm at each frame normalized by the average rMSE for the comparison algorithm across the whole video clip. The median improvement is calculated in the same manner. 

The recovery results for this video database show consistent performance improvements for RWL1-DF when evaluated over all video sequences in this database.  Additionally, we note that some video sequences were significantly richer in texture and motion than others, resulting in a more challenging recovery task. We identified 13 such video clips that were especially challenging (i.e, those where the average rMSE reconstruction for BPDN is over 1\%).  For these clips we plot the mean and median percent improvement in Figure~\ref{fig:ImpBar}.  RWL1-DF shows very significant improvements within these video sequences, indicating that RWL1-DF is especially beneficial in challenging recovery scenarios.  Additionally, an adaptation of a similar RWL1 method to a spatial filtering setting has shown similar performance improvements in the most challenging data examples~\cite{charles2014rwl1sf}.

\section{conclusions}

In this work we explore two causal dynamic filtering schemes for time-varying sparse signals that are based on \lo minimization approaches. We examine both a simple method (BPDN-DF) that admits a theoretical analysis as well present a novel algorithm (RWL1-DF) that achieves state-of-the-art performance by leveraging second-order statistics in a similar manner as Kalman filtering.  We conclude from these results that \lo minimization methods have significant potential to be extended from static sparse signal estimation to time-varying data streaming settings.   While not discussed in detail, the approaches in this paper can be solved efficiently through simple approaches such as an Iterative Shrinkage and Thresholding Algorithm (ISTA). Any future advances in numerical algorithms for \lo optimization will allow for faster and larger scale implementations of these approaches.

There are a number of potential modifications and improvements to the described approaches that will be pursued in future work.  First, while classic dynamic filtering approaches assume known dynamics, it would be natural to estimate or learn complex dynamics as part of the proposed model structure.  Second, this paper follows the convention of classic dynamic filtering work by assuming that the dynamics model operates in the state space (i.e., predicting the variances from the previous state estimates).  In scenarios where the dynamics operations are not derived from physical principles and may be learned from data, RWL1-DF may be more effective when the dynamics model predicts the variances from the previous variance estimates.


\section{Appendix}

\subsection{Temporal bound derivation}

\label{app:bpdndfproof}

This appendix proves the result in Theorem~\ref{thm:BPDNDFbasic}. While the results do not depend on the numerical algorithm used to solve the optimization, the proof technique will leverage ISTA in the analysis~\cite{bredies2008linear,DAU:2004,aureleISTA}. 
This proof will follow by using the ISTA update equation at each time-step $n$ and algorithmic iteration $l$, along with the RIP condition on $\Phi$ and Lemma~\ref{lem:ubound} (presented below), to obtain a recurrence relation for the absolute error between the current signal estimate. 
We then solve the recurrence relation for the analytical upper bound for each $l,n$ independent of the previous ISTA iteration $l+1$ (but still dependent on the previous time-step $n-1$), allowing us to take $l\rightarrow\infty$ and obtaining a recursive steady state bound in terms of only $n$.
We then solve the new recursion equation in terms of $n$, yielding an upper bound for all $n$ depending only on the initial condition and the properties of the dynamics and measurement functions. 

We start by noting that ISTA for BPDN-DF can can be written as
\begin{eqnarray}
	\istavect{n}^{l+1} & = & \coefvtest{n}^l + \frac{\ssize}{1+\kthresh}\basemat^T\left(\mematt{n}^T\left(\insigt{n} - \mematt{n}\basemat\coefvtest{n}^l\right)\right. \nonumber \\
	& & \qquad\left. + \kthresh\left(f(\basemat\coefvtest{n-1}) - \basemat\coefvtest{n}^{l}\right) \right) \nonumber \\
	\coefvest^{l+1} & = & T_{\gamma}(\istavect{n}^{l+1}),  \label{eqn:ista}
\end{eqnarray}
where \istavec is an un-thresholded state, $T_{\gamma}(\cdot)$ is a soft-thresholding function with threshold parameter $\gamma$, and \ssize is the algorithm's step size. 
The proof continues by demonstrating that error between the estimate and the true signal contracts at each algorithmic step. 
Using the fact that the converged solution of ISTA minimizes the corresponding cost function~\cite{bredies2008linear,DAU:2004}, the error bound on the ISTA solution is then also a bound for the overall BPDN-DF minimizer.
For ease of notation, we omit temporal subscripts, assuming that all variables, unless otherwise stated, have temporal subscript $n$. Additionally, we define the previous (steady state) estimate as $\coefvtest{n} = \coefprev$ and refer to $\coeftrue_n$ as the true coefficients at time $n$. We define two subsets of the sparse vector \coefvec: $J$ and $J^\prime$. The set $J = J[l+1]$ is the union of the current set of active coefficients in the ISTA algorithm $\Gamma[l]$, the $q$ largest elements of the vector $\istavec$ $\Delta[l]$, and the true active set $\Gamma_{\dagger}$. In~\cite{aureleISTA} it is shown that $|J| = |\Delta[l+1]\cup\Gamma[l]\cup\Gamma_\dagger| \leq \sparsity + 2q$. Similarly we define $J^\prime = J[l+2] = \Delta[l+1]\cup\Gamma[l+1]\cup\Gamma_\dagger$. To start, we require the following lemma which bounds the norm of $\istavec$ at each algorithmic step and at each $n$:

\begin{lemma}
    \label{lem:ubound}
    Suppose that the same conditions as in Theorem~\ref{thm:BPDNDFbasic} hold. Additionally, assume that $\norm{\coefvect{n}^\dagger}_2^2 \leq b$ for all $n$. The vector $\istavect{J}^l$ (the ISTA variables $\istavec^l$ restricted to the support subset $J$) at each algorithmic iteration $l$ obtained via ISTA (iterating Equation~\eqref{eqn:ista}) with a step size of $\stepsz$ satisfies  
    \begin{gather}
	    \|\istavect{J}^l\|_2 \leq \left(|\ssize - 1| + \frac{\ssize\ripd}{1+\kthresh}\right)\frac{\gamma}{1+\kthresh}\sqrt{q} + \left(\ssize + \ssize\frac{\kthresh f^{\ast}  + \ripd}{\kthresh + 1}\right)b  \nonumber \\
	    +  \ssize\frac{\sqrt{1+\delta}}{1+\kthresh}\errmax + \frac{\ssize\kthresh}{1+\kthresh}\inmax, \nonumber 
    \end{gather}
    And with the restriction
    \begin{gather}
	\ssizek(\kthresh + \kthresh\dlip + 1 +\ripd)b + \ssizek\sqrt{1+\ripd}\errmax + \ssizek\kthresh\inmax \nonumber \\
	\qquad\qquad\qquad  \leq \left(1- \frac{|\ssizek(1+\kthresh) - 1| + \ssizek\ripd)}{1+\kthresh}\right)\gamma\sqrt{q}, \nonumber
    \end{gather}
    then $\norm{\istavect{J}^l}_2$ is simply bounded by $\norm{\istavect{J}^l}_2  \leq  \gamma\sqrt{q}$
\end{lemma}

We start by writing the norm using the definition of \istavec in the ISTA algorithm:
\begin{gather}
	\norm{\istavect{J}^l}_2  =  \left\|\coefvect{J}^l + \frac{\ssize}{1+\kthresh}(\memat\basemat)_J^T\left(\insig - \memat\basemat\coefvec^l\right)\right.  \nonumber \\
	\qquad\qquad \left. + \frac{\ssize\kthresh}{1+\kthresh}\basemat_J^T\left( f(\basemat\coefprev) - \basemat\coefvec^l\right) \right\|_2. \nonumber
\end{gather}
Using the fact that $\insig = \memat\basemat\coeftrue + \meerrv$,
\begin{gather}
	\norm{\istavect{J}^l}_2  =  \left\| \coefvect{J}^l + \frac{\ssize}{1+\kthresh}(\mematt\basemat)_J^T\memat\basemat\left(\coeftrue - \coefvec^l\right) \right. \nonumber \\
	\left. + \frac{\ssize}{1+\kthresh}\basemat^T\mematt{J}^T\meerrv + \frac{\ssize\kthresh}{1+\kthresh}\basemat_J^T\left( f(\basemat\coefprev) - \basemat\coefvec^l \right)   \right\|_2. \label{eqn:ubound}
\end{gather}

To properly reduce the portion of this expression depending on the dynamics $f(\cdot)$, we note that since the dynamics satisfies $
\basemat\coeftrue_n = f(\basemat\coeftrue_{n-1}) + \innovvt{n}$, 
\begin{eqnarray}
	f(\basemat\coefprev) - \basemat\coefvec^l & = & f(\basemat\coefprev) - \basemat\coeftrue + \basemat\coeftrue - \basemat\coefvec^l \nonumber \\
	& = & f(\basemat\coefprev) - f(\basemat\coeftrue_{n-1}) - \innovv \nonumber \\
	& & \qquad\qquad + \basemat\left(\coeftrue - \coefvec^l\right). \nonumber
\end{eqnarray}
Setting $\ssizek = \ssize/(1+\kthresh)$ and collecting similar terms,
\begin{eqnarray}
	\norm{\istavect{J}^l}_2 & = & \left\| \coefvect{J}^l + \ssizek\left( (\memat\basemat)_J^T\memat\basemat + \kthresh\bm{I}_J \right)\left(\coeftrue - \coefvec^l \right) \right. \nonumber \\
	& &  + \ssizek\basemat^T\left(\mematt{J}^T\meerrv - \kthresh\innovv\right) \nonumber \\
	& & \left. + \ssizek\kthresh\basemat_J^T\left( f(\basemat\coefprev) - f(\basemat\coeftrue_{n-1})\right)   \right\|_2 \nonumber \\
	& \leq & \norm{\coefvect{J}^l + \ssizek((\memat\basemat)_J^T\memat\basemat + \kthresh\bm{I}_J)\left(\coeftrue - \coefvec^l\right)}_2  \nonumber \\
	& & \qquad + \ssizek\sqrt{1+\ripd}\norm{\meerrv}_2  + \ssizek\kthresh\norm{\innovv}_2 \nonumber \\ 
	& & \qquad + \ssizek\kthresh\norm{\basemat_J^T\left( f(\basemat\coefprev) - f(\basemat\coeftrue_{n-1})\right) }_2 \nonumber \\
	& \leq & \norm{ \coefvect{J}^l + \ssizek( (\memat\basemat)_J^T\memat\basemat + \kthresh\bm{I}_J)\left(\coeftrue - \coefvec^l\right) }_2  \nonumber \\
	& & + \ssizek\sqrt{1+\ripd}\norm{\meerrv}_2  + \ssizek\kthresh\norm{\innovv}_2 + \ssizek\kthresh\dlip\norm{\bm{e}_{n-1} }_2 \nonumber \\
	& \leq & \norm{ (\ssizek(\memat\basemat)_J^T\memat\basemat - (1 - \kthresh\ssizek)\bm{I}_J)\coefvect{J}^l}_2 \nonumber \\
	& & + \norm{\ssizek(\mematt{J}^T\memat + \kthresh\bm{I}_J)\basemat\coeftrue}_2 + \ssizek\sqrt{1+\ripd}\norm{\meerrv}_2 \nonumber \\
	& & \qquad\qquad + \ssizek\kthresh\norm{\innovv}_2 + \ssizek\kthresh\dlip\norm{\bm{e}_{n-1}} _2.    \label{eqn:interbound1}
 \end{eqnarray}
 where the first and third inequalities follow from the triangle inequality, the fact that $\norm{\basemat} \leq 1$ and the RIP of \memat, and the second inequality follows from the smoothness of $f(\cdot)$. 

We now use the Cauchy-Schwartz inequality with the fact that, since \memat satisfies RIP$(|J|,\ripd)$ with respect to \basemat, 
\begin{gather}
	\norm{\alpha(\basemat\memat)_J^T(\memat\basemat)_J + \beta\bm{I}_J}_2 \leq |\alpha + \beta| + \alpha\ripd, \label{eqn:eigineq}
\end{gather}
for any constants $\alpha$, $\beta$, to bound the first two terms of Equation~\ref{eqn:interbound1} as  
\begin{gather}
	\norm{\istavect{J}^l}_2  \leq (|\ssizek + \kthresh\ssizek - 1| + \ssizek\ripd )\norm{\coefvec^l}_2 + \ssizek(\kthresh + 1 + \ripd)\norm{\coeftrue}_2 \nonumber \\
	\qquad \qquad + \ssizek\sqrt{1+\ripd}\norm{\meerrv}_2 + \ssizek\kthresh\norm{\innovv}_2 + \ssizek\kthresh\dlip\norm{\bm{e}_{n-1}}_2, \nonumber 
\end{gather}

Simplifying, we obtain
\begin{eqnarray}
	\norm{\istavect{J}^l}_2 & \leq & \left(|\ssizek(1+\kthresh) - 1| + \ssizek\ripd\right)\widetilde{\gamma}\sqrt{q} + \ssizek(\kthresh + 1 + \ripd)b  \nonumber \\
	& & \qquad  + \ssizek\sqrt{1+\ripd}\errmax + \ssizek\kthresh\inmax + \ssizek\kthresh\dlip\norm{\bm{e}_{n-1}}_2, \nonumber
\end{eqnarray}
where $b$ is the maximum energy of $\coeftrue$ (i.e. $\norm{\coeftrue}_2 \leq b$) and $\widetilde{\gamma} = \gamma/(1+\kappa)$. 

Ideally, this Lemma should be independent of the previous estimation error norm $\norm{\bm{e}_{n-1}}_2$ for the Lemma to hold for all $n$. If we initialize the estimate with the zero vector, clearly $\norm{\bm{e}_0}_2 = \norm{\coeftrue_0}_2 \leq b$. Thus, setting $\norm{\bm{e}_{n-1}}_2 < b$ results in the Lemma statement and it remains to ensure that $\norm{\bm{e}_n}_2 \leq b$ for all $n$ through later parameter restrictions. 

\subsection{First Recursion}
With Lemma~\ref{lem:ubound}, we seek a recursive expression for the estimation error at algorithmic iteration $l+1$,
\begin{eqnarray}
	\norm{\bm{e}^{l+1}}_2 & = & \norm{\coefvec^{l+1} - \coeftrue}_2 \nonumber \\
	& = & \norm{\coefvec^{l+1} - \istavect{J^\prime}^{l+1} + \istavect{J^\prime}^{l+1} - \coeftrue}_2 \nonumber \\
	& \leq & \norm{\coefvec^{l+1} - \istavect{J^\prime}^{l+1}}_2 + \norm{\istavect{J^\prime}^{l+1} - \coeftrue}_2 \nonumber \\
	& \leq & \norm{\istavect{J^\prime}^{l+1}}_2 + \left\|\coefvect{J^\prime} - \coeftrue \right. \nonumber \\
	& & + \ssizek\basemat^T\left(\mematt{J}^T\memat\basemat(\coeftrue - \coefvect{J^\prime}^l) \right. \nonumber \\
	& & \left.\left. + \mematt{J^\prime}\meerrv + \kthresh(f(\basemat\coefprev) - \basemat\coefvec^l) \right) \right\|_2 \nonumber \\
	& = & \norm{\istavect{J^\prime}^{l+1}}_2  + \left\|(\ssizek\basemat^T\mematt{J}^T\memat\basemat - \bm{I}_{J^\prime})(\coeftrue-\coefvect{J^\prime}^l) \right. \nonumber \\
	& & \left. + \ssizek\basemat^T\mematt{J^\prime}\meerrv + \ssizek\kthresh\basemat^T(f(\basemat\coefprev) - \basemat\coefvec^l)  \right\|_2 \nonumber \\
	& \leq & \gamma\sqrt{q} + \ssizek\kthresh\dlip\norm{\bm{e}_{n-1}}_2  + \ssizek\sqrt{1+\ripd}\errmax + \ssizek\kthresh\inmax  \nonumber \\
	& & \qquad + (|\ssizek + \ssizek\kthresh - 1| + \ssizek\ripd)\|\bm{e}^l\|_2, \nonumber 
\end{eqnarray}
where the first inequality follows from the triangle inequality, the second inequality follows from the nature of the thresholding function and the definition of \istavec, the third inequality follows from Lemma~\ref{lem:ubound} for the first term and a similar set of steps as used to prove Lemma~\ref{lem:ubound} used on the second term. This recursive formula can now be solved for $\norm{\bm{e}^l}_2$,
\begin{eqnarray}
	\norm{\bm{e}[l]}_2  & \leq & (|\ssize - 1| + \ssizek\ripd)^l\left(\norm{\bm{e}^0}_2  \right. \nonumber \\
	& & \left. - \frac{\gamma\sqrt{q} + \ssizek\kthresh\dlip\norm{\bm{e}_{n-1}}_2  + \ssizek\sqrt{1+\ripd}\errmax + \ssizek\kthresh\inmax}{1 - |\ssize - 1| - \ssizek\ripd}\right) \nonumber \\
	& & + \frac{\gamma\sqrt{q} + \ssizek\kthresh\dlip\norm{\bm{e}_{n-1}}_2  + \ssizek\sqrt{1+\ripd}\errmax + \ssizek\kthresh\inmax}{1 - |\ssize - 1| - \ssizek\ripd}. \nonumber 
\end{eqnarray}
As $l\rightarrow \infty$, ISTA converges to a solution with a bounded error of 
\begin{eqnarray}
	\norm{\bm{e}_n}_2 & = & \norm{\bm{e}_n^{\infty}}_2 \nonumber \\
	& \leq & \frac{\gamma\sqrt{q} + \ssizek\kthresh\dlip\norm{\bm{e}_{n-1}}_2  + \ssizek\sqrt{1+\ripd}\errmax + \ssizek\kthresh\inmax}{1 - |\ssize - 1| - \ssizek\ripd}. \nonumber
\end{eqnarray}

\subsection{Second Recursion}
The steady state error of ISTA can now be treated as a new recursive equation on $\norm{\bm{e}_n}_2$ in terms of the signal time-step $n$. This expression can be solved for $\norm{\bm{e}_n}_2$
\begin{eqnarray}
	\norm{\bm{e}_n}_2  & \leq & \left(\frac{\ssizek\kthresh\dlip}{1 - |\ssize - 1| - \ssizek\ripd}\right)^n\left(\norm{\bm{e}_0}_2 \right. \nonumber \\
	& & \left. - \frac{\gamma\sqrt{q} + \ssizek\sqrt{1+\ripd}\errmax + \ssizek\kthresh\inmax}{1 - |\ssize - 1| - \ssizek\ripd - \ssizek\kthresh\dlip} \right) \nonumber \\
	 & & \qquad\qquad + \frac{\gamma\sqrt{q} + \ssizek\sqrt{1+\ripd}\errmax + \ssizek\kthresh\inmax}{1 - |\ssize - 1| - \ssizek\ripd - \ssizek\kthresh\dlip}, \nonumber 
\end{eqnarray}
which yields a bound for the error at every iteration $n$.

The only remaining task is to ensure $\norm{\bm{e}_n}_2 \leq b$ for all $n$ so that Lemma~\ref{lem:ubound} holds. If $b$ bounds $\norm{\bm{e}_n}_2$ at each $n$, then
\begin{gather}
	\beta^n(\norm{\bm{e}_0}_2 - \frac{\alpha}{1-\beta}) + \frac{\alpha}{1-\beta} \leq \beta^n(b - \frac{\alpha}{1-\beta}) + \frac{\alpha}{1-\beta} \leq b \nonumber 
\end{gather}
holds.
Simplifying,
\begin{gather}
	\frac{\alpha(1-\beta^n)}{1-\beta} \leq (1-\beta^n)b \rightarrow \frac{\alpha}{1-\beta} \leq b, \nonumber
\end{gather}
or in terms of the model and system parameters,
\begin{gather}
\frac{(1+\kthresh)\gamma\sqrt{q} +  \ssize\sqrt{1+\ripd}\errmax + \ssize\kthresh\inmax}{(1 + \kthresh)(1 - |\ssize - 1|) - \ssize\ripd - \ssize\kthresh\dlip} \leq b. \nonumber
\end{gather}
Setting $\ssize = (1+\ripd)^{-1} \leq 2(1+\ripd)^{-1}$ completes the proof.


\bibliographystyle{IEEEtran}

\bibliography{IEEEabrv,RWL1,MPbib,CSbib,BPDNBIB,Dynamicbib,Kalmanbib,AudioBib,Adam_Bibliography2}

\end{document}